\documentclass[10pt,twoside]{article}
\usepackage{amsmath}
\usepackage{Latex-document}
\usepackage{graphicx}
\usepackage{amsfonts}
\usepackage{amssymb}
\newtheorem{theorem}{Theorem}
\newtheorem{conjecture}[theorem]{Conjecture}

\newtheorem{definition}[theorem]{Definition}

\title{\textbf{Automorphism Groups of Saturated \vskip-2mm Structures; A Review\vskip6mm}}
\author{D. Lascar\vspace*{-0.5cm}\thanks{CNRS, Universit\'e Denis Diderot Paris 7, 2
Place Jussieu, UFR de math\'ematiques, case 7012, 75521 Paris
Cedex 05, France. E-mail: lascar@logique.jussieu.fr}}
\date{\vspace{-8mm}}

\begin{document}

\maketitle

\thispagestyle{first} \setcounter{page}{25}

\begin{abstract}

\vskip 3mm

We will review the main results concerning the automorphism groups of saturated structures which were obtained
during the two last decades. The main themes are: the small index property in the countable and uncountable cases;
the possibility of recovering a structure or a significant part of it from its automorphism group; the subgroup of
strong automorphisms.

\vskip 4.5mm

\noindent\textbf{2000 Mathematics Subject Classification:} 03C50, 20B27.

\noindent {\bf Keywords and Phrases:} Automorphism groups, Small index property, Strong automorphisms.
\end{abstract}

\vskip 12mm

\section{Introduction}

\vskip-5mm \hspace{5mm}

Saturated models play an important role in model theory. In fact, when
studying the model theory of a complete theory $T$, one may work in a large
saturated model of $T$ with its definable sets, and forget everything else
about $T$. This large saturated structure is sometimes called the ``universal
domain'', sometimes the ``monster model''.

A significant work has been done the last twenty years on the automorphism
groups of saturated models. It is this work that I want to review here.
There is a central question that I will use as a ``main theme'' to organize
the paper: what information about $M$ and its theory are contained in its
group of automorphisms? In the best case, $M$ itself is ``encoded'' in some
way in this group; recovering $M$ from it is known as ``the reconstruction
problem''. A possible answer to this problem is a theorem of the form: If
$M_{1}$ and $M_{2}$ are structures in a given class with isomorphic
automorphism groups, then $M_{1}$ and $M_{2}$ are isomorphic.

Throughout this paper, $T$ is supposed to be a countable complete theory. The countability of $T$
is by no means an essential hypothesis. Its purpose is only to make the exposition
smoother, and most of the results generalize without difficulty to
uncountable theories. We will denote by $Aut(M)$ the group of automorphisms of
the structure $M$, and if $A$ is a subset of $M$, $Aut_{A}(M)$ will be the
pointwise stabilizer of $A$:
\[
Aut_{A}(M)=\left\{  f\in Aut(M)\;;\;\forall a\in
A\,f(a)=a\right\}.
\]

When we say ``definable", we mean ``definable without parameters".

\section{The countable case}

\vskip-5mm \hspace{5mm}

As a preliminary remark, let us say that the automorphism group of a
saturated model is always very rich: if $M$ has cardinality $\lambda$, then
its automorphism group has cardinality $2^{\lambda}$.

I do not know who was the first to introduce the small index property. As we
will see, it is crucial in the subject.

\begin{definition}
Let $M$ be a countable structure. We say that $M$ (or $Aut(M)$) has the small
index property if for any subgroup $H$ of $Aut(M)$ of index less than
$2^{\aleph_{0}}$, there exists a finite set $A\subset M$ such that
$Aut_{A}(M)\subseteq H$.
\end{definition}

Remark that the converse is true: any subgroup containing a subgroup of the
form $Aut_{A}(M)$ where $A$ is finite, has a countable index in $Aut(M)$%
. Moreover, the subgroups containing a subgroup of the form $Aut_{A}(M)$ are
precisely the open neighborhoods of the identity for the pointwise convergence
topology. In other words, the small index property allows us to recover the
topological structure of $Aut(M)$ from its pure group structure.

The small index property has been proved for a number of countable saturated structures:

\begin{enumerate}
\item  The infinite set without additional structure \cite{semmes}, \cite{dixon}.

\item  The linear densely ordered sets \cite{truss}.

\item  The vector spaces over a finite field \cite{evan86}.

\item  The random graph \cite{HHLS}.

\item  Various other classes of graphs \cite{herwig98}.

\item  Generic relational structures \cite{herwig95}.

\item $\omega$-categorical $\omega$-stable structures \cite{HHLS}.
\end{enumerate}

The small index property has also been proved for some countable structures
which are not saturated: for the free group with
$\omega$-generators (\cite{briantevans}), for arithmetically saturated models of
arithmetic (\cite{lascar94}).

There are examples of countable saturated structures which fail to have the
small index property. The simplest may be an algebraically closed field of
characteristic 0 of infinite countable transcendence degree: Let
$\mathbb{\bar{Q}}$ be the algebraic closure of the field of rational
numbers.\ There is an obvious homomorphism $\varphi$ from $Aut(M)$ onto
$Aut(\mathbb{\bar{Q}})$ (the restriction map). Now, it is well known that
there is a subgroup $H$ of $Aut(\mathbb{\bar{Q}})$ of countable index (in fact
of finite index) which is not closed for the Krull topology, which is nothing
else that the pointwise convergence topology.\ Then $\varphi^{-1}(H)$ is not
open, but of finite index in $Aut(M)$.

As we will see later, the small index property is particularly relevant for
$\omega$-categorical structures.\ Evans and Hewitt have produced an example of
such a structure without the small index property (\cite{evanshewitt}).

With the pointwise convergence topology, $Aut(M)$ is a topological polish
group.\ So, we may use the powerful tools of descriptive set theory. In
many cases (for example for structures 1-6 above), it can be shown that there
is a (necessarily unique) conjugacy class which is generic, that is, is the
countable intersection of dense open subsets. The elements of this class are
called generic automorphisms, and they play an important role in the proof of
the small index property.

Another possible nice property of these automorphism groups which is sometimes
obtained as a bonus of the proof of the small index property, is the fact that
its cofinality is not countable, that is, $Aut(M)$ is not the union of a
countable chain of proper subgroups. This is proved in particular for the full
permutation group of a countable set (\cite{macphersonneumann}), for the random
graph and for $\omega$-categorical $\omega$-stable structures (\cite{HHLS}).

I would like to mention here the work of Rubin (\cite{rubin}). He has shown how to
reconstruct a certain number of structures from their automorphism group
using a somewhat different method.\ His methods apply essentially to
``combinatorial structures'' such as the random graph, the universal
homogeneous poset, the generic tournament (a structure for which the small
index property is not known), etc.

\section{Subgroups and imaginary elements}

\vskip-5mm \hspace{5mm}

Recall that an imaginary element of $M$ is a class of a tuple of $M^n$ modulo a definable
equivalence relation on $M^n$. For instance, if $G$ is a group and $H$ a definable subgroup
of $G^n$, then any coset of $H$ in $G^n$ is an imaginary element. When we add all
these imaginary elements to a saturated structure
$M$, we obtain the structure $M^{eq}$, and we can consider $M^{eq}$ as a
saturated structure (in a larger language).

It is clear that $M$ and $M^{eq}$ have canonically the same
automorphism group: every automorphism of $M$ extends uniquely to an automorphism of $M^{eq}$.
 This shows a limitation to the reconstruction problem: If
$M$ and $N$ are two structures which are such that ``$M^{eq}$ and $N^{eq}$ are
isomorphic'', then $Aut(M)$ and $Aut(N)$ are isomorphic via a bicontinuous
isomorphism. The condition ``$M^{eq}$ and $N^{eq}$ are isomorphic'' may seem
weird, but in fact, it is natural. Roughly speaking, it means that $M$ can be
interpreted in $N$, and conversely (a little more in fact, see
\cite{Ahlbrandtziegler} for more details). In this case, we say that $M$ and
$N$ are bi-interpretable.

Consider now the case of an $\omega$-categorical structure $M$. It is not
difficult to see that any open subgroup of $Aut(M)$ is the stabilizer
$Aut_{\alpha}(M)$ of an imaginary element $\alpha$. Moreover, $Aut(M)$ acts by
conjugation on the set of its open subgroups, and this action is (almost)
isomorphic to the action of $Aut(M)$ on $M^{eq}$ (almost because two different
imaginary elements $\alpha$ and $\beta$ may have the same stabilizer). So,
from the topological group $Aut(M)$ we can (almost) reconstruct its action on
$M^{eq}$. We can do better:

\begin{theorem}
\emph{\cite{Ahlbrandtziegler} }Assume that $M$ and $N$ are countable $\omega
$-categorical structures. Then the following two conditions are equivalent:

\begin{enumerate}
\item  there is a bicontinuous isomorphism from $Aut(M)$ onto $Aut(N);$

\item $M$ and $N$ are bi-interpretable.
\end{enumerate}
\end{theorem}

In fact, these conditions are also equivalent to: there exists a continuous
isomorphism from $Aut(M)$ onto $Aut(N)$ (see \cite{lascar91}). Thus, if one of
the structure $M$ or $N$ has the small index property and $Aut(M)$ is
isomorphic to $Aut(N)$ (as pure groups), then $M$ and $N$ are bi-interpretable.

Now, if $M$ is not necessarily $\omega$-categorical (but still saturated), the
situation is a bit more complicated. We need to introduce new elements.

\begin{definition}
\noindent 1. An ultra-imaginary element of $M$ is a class modulo $E$, where $E$
is an equivalence relation on $M^{n}$ ($n\leq\omega)$ which is invariant under
the action of $Aut(M)$. An ultra-imaginary element is finitary if $n<\omega$.

\noindent 2. A hyperimaginary element of $M$ is a class modulo $E$, where $E$
is an equivalence relation on $M^{n}$ ($n\leq\omega)$ which is defined by a
(possibly infinite) conjunction of first order formulas.
\end{definition}

An imaginary element is hyperimaginary, and a hyperimaginary element is
ultra-imaginary. A hyperimaginary element is a class modulo an equivalence
relation $E$ defined by a formula of the form
\[
\bigwedge_{i\in I}\varphi_{i}%
\]
where the $\varphi_{i}$ are first-order formulas (without parameters) and whose
free variables are among the $x_{k}$ for $k<n$. An ultra-imaginary element is a
class modulo an equivalence relation $E$ defined by a formula of the form
\[
\bigvee_{j\in J}\bigwedge_{i\in I}\varphi_{i,j}%
\]
where the $\varphi_{i,j}$ are first order-formulas (without parameters) and
whose free variables are among the $x_{k}$ for $k<n$.

If $M$ is a countable saturated structure, the stabilizer of a finitary
ultra-imaginary element is clearly an open subgroup, and it is not difficult to
see that if $H$ an open subgroup of $Aut(M)$, then there exists a finitary
ultra-imaginary element $\alpha$ such that $H$ is the stabilizer of $\alpha$.
In the $\omega$-categorical case, any finitary ultra-imaginary is in fact
imaginary, and this explain why this case is so simple.

In some cases, for example for $\omega$-stable theories (see \cite{lascar96}),
it is possible to characterize, among all open subgroups, those which are of
the form $Aut_{\alpha}(M)$ with $\alpha$ imaginary. Something similar has been
done for countable arithmetically saturated models of arithmetic in
\cite{kotlarskikaye}, and in \cite{kossakschmerl}, it is proved that if two such models have
isomorphic automorphism groups, then they are isomorphic.

\section{Strong automorphisms}

\vskip-5mm \hspace{5mm}

It is now time to introduce the group of strong automorphisms.

\begin{definition}
\emph{\cite{lascar82}} The group of strong automorphisms of $M$ is the group
generated by the set
\[
\bigcup\left\{  Aut_{N}(M)\;;\;N\prec M\right\}
\]
and is denoted $Autf(M)$.
\end{definition}

It is easy to see that $Autf(M)$ is a normal subgroup of $Aut(M)$. Its index
is at most $2^{\aleph_{0}}$.\ Moreover, the quotient group $Aut(M)/Autf(M)$
depend only on $T$ : if $M$ and $M^{\prime}$ are two saturated models,
$M \prec M^{\prime}$, then there is a natural isomorphism from $Aut(M)/Autf(M)$ onto
$Aut(M^{\prime})/Autf(M^{\prime})$. $Aut(M)/Autf(M)$ will be denoted $Gal(T) $
(of course, $Gal$ stands for Galois). For example, if $T$ is the theory of
algebraically closed fields of characteristic 0, $Autf(M)=Aut_{\mathbb{\bar
{Q}}}(M)$ and $Gal(T)$ is (isomorphic to) the group of automorphisms of
$\mathbb{\bar{Q}}$.

In fact this interpretation is general. Assume first that $M$ is of
cardinality bigger than $2^{\aleph_{0}}$. Let $\alpha$ be an ultra-imaginary
element of $M$. It can be shown that the following conditions are equivalent:

\begin{enumerate}
\item $card\{f(\alpha)\;;\;f\in Aut(M)\}<card(M)$;

\item $card\{f(\alpha)\;;\;f\in Aut(M)\}\leq2^{\aleph_{0}}$; \newline An
equivalence relation is bounded if it has at most $2^{\aleph_{0}}$ classes
(equivalently less than $card(M)$ classes). The above conditions are also
equivalent to:

\item $\alpha$ (as a set) is the class modulo an invariant bounded equivalence relation.
\end{enumerate}

If these conditions are satisfied, we say that $\alpha$ is bounded. It should
be remarked that an imaginary element is bounded if and only if it is algebraic,
if and only if its orbit is finite.

We will denote by $Bdd(M)$ the set of bounded ultra-imaginary elements of $M$.
This set does not really depend on $M$ (as soon as its cardinality is big
enough) but only on its theory: any invariant bounded equivalence relation has
a representative in any uncountable saturated model.\ We will allow ourself to
write $Bdd(T)$ when convenient. Moreover $Autf(M)$ is exactly the pointwise
stabilizer of $Bdd(M)$ so that $Gal(T)$ can be identified with the group of
elementary permutations of $Bdd(T)$.\

With some care, we can generalize this interpretation to models of small
cardinality: for example, assume $M$ countable, and let $M^{\prime}$ be
a large saturated extension of $M$. Then any automorphism $f$ of $M$ extends to
an automorphism of $M^{\prime}$, and if $f_{1}$ and $f_{2}$ are two such
extensions, then their action on $Bdd(M^{\prime})$ are equal; $Autf(M)$ is
exactly the set of automorphisms whose extensions to $M^{\prime}$ act
trivially on $Bdd(M^{\prime})$.

In any case, $Aut(M)$ leaves fixed the set of bounded imaginary elements and
the set of bounded hyperimaginary elements. In some cases (for example for
algebraically closed fields), $Gal(T)$ acts faithfully on the set of bounded
imaginary elements. It is the case if $T$ is stable (\cite{lascar82}). It is
not known if it is always true for simple theories, but it is true for the
so-called low simple theories (\cite{buechler}) and in particular for
supersimple theories. For simple theories, $Gal(T)$ acts faithfully on the set
of bounded hyperimaginary elements (\cite{kim}). In \cite{casanovasandall}
there is an example of a theory where the action of $Gal(T) $ on the set of
hyperimaginary elements is not faithful.

There is a natural topology on $Gal(T)$ (see \cite{lascarpillay} for details).
It can be defined in two different ways.\

My favorite one is via the ultraproduct construction.\ Let $(\gamma
_{i}\;;\;i\in I)$ be a family of elements of $Gal(T)$ and $\mathcal{U}$ an
ultrafilter on $I$.\ Choose a saturated model $M$ and, for each $i\in I$ an
automorphism $f_{i}\in Aut(M)$ lifting $\gamma_{i}$. Consider the
ultrapower $M^{\prime}=$  $\prod_{i\in\mathcal{U}}M$. We can define the
automorphism $\prod_{i\in\mathcal{U}}f_{i}$ on $M^{\prime}$.
This automorphism acts on $Bdd(M^{\prime})=Bdd(T)$, so defines an element of
$Gal(T)$, say $\beta$. This element $\beta$ should be considered as a limit of
the family $(\gamma_{i}\;;\;i\in I)$ along $\mathcal{U}$.\ A subset $X$ of
$Gal(T)$ is closed for the topology we are defining if it is closed for this
limit operation. You should be aware that the element $\beta$ may depend on
the choices of the $f_{i}$'s, because the topology we are defining is not
Hausdorff in general.

The other way to define a topological structure on $Gal(T)$ is to define a
topology on $Bdd(T)$. If, as it is the case when $T$ is stable, $Gal(T)$ can
be identified with a group of permutation on the set of imaginary elements,
then we just endow $Gal(T)$ with the pointwise convergence topology
(that is we consider the set of imaginary elements with the discrete
topology). Otherwise, it is more complicated, and here is what
should be done in general:

For each $n\leq\omega$ and $E$ invariant bounded equivalence relation on
$M^{n}$, consider the canonical mapping $\varphi_{E}$ from $M^{n}$ onto
$M^{n}/E$. By definition, a subset $X$ of  $M^{n}/E$ is closed if and only if
$\varphi_{E}^{-1}(X)$ is the intersection of a family of subsets definable with
parameters. $Gal(T)$ acts on $M^{n}/E$ and the topology on $Gal(T)$ is defined
as the coarsest topology which makes all these actions (with various $n$ and
$E$) continuous.

Now, we can prove:

\begin{theorem} 1. $Gal(T)$ is a topological compact group.

\begin{enumerate}
\item[2.]  It is Hausdorff if and only if it acts faithfully on the set of bounded
hyperimaginary elements, if and only if it acts faithfully on the set of finitary
bounded hyperimaginary elements.

\item[3.]  It is profinite if and only if it acts faithfully on the set of
bounded imaginary elements.
\end{enumerate}
\end{theorem}

There is a Galois correspondence between the subgroups of $Gal(T)$ and the bounded
ultra-imaginary elements: every subgroup of $Gal(T)$ is the stabilizer of an ultra-imaginary element. The hyperimaginary elements correspond to the closed
subgroups and the imaginary elements correspond to the clopen subgroups of $Gal(T)$.

Let $H_{0}$ be the topological closure of the identity. Then $H_{0}$ is a
normal subgroup of $Gal(T)$. If we consider $Gal(T)$ as a permutation group on
$Bdd(T)$, $H_{0}$ is exactly the pointwise stabilizer of the set of bounded
hyperimaginary elements. So, if we set $Gal_{0}(T)=Gal(T)/H_{0}$, $Gal_{0}(T)$
acts faithfully on the set of bounded hyperimaginary elements. As a quotient
group, $Gal_{0}(T)$ is canonically endowed with a topology. This way, we get a
compact Hausdorff group.\

Recently, L.\ Newelski (\cite{newelski}) has proved that $H_{0}$ is either
trivial or of cardinality $2^{\aleph_{0}}$.

I would like to conclude this section by a conjecture. In all the known examples of
countable saturated structures where the small index property is false, there
is a non open subgroup of $Gal(T)$ of countable index (and, its preimage
by the canonical homomorphism from $Aut(M)$ onto $Gal(T)$ is a non open
subgroup of $Aut(M)$ of countable index). If $A$ has a cardinality strictly
less than $card(M)$, define $Autf_{A}(M)$ as the subgroup of $Aut_{A}(M)$
generated by
\[
\bigcup\left\{  Aut_{N}(M)\;;\;A\subseteq N\prec M\right\}.
\]
The following conjecture is open, even in the $\omega$-categorical case:

\begin{conjecture}
\label{conjecture}Assume that $M$ is a countable saturated structure and let
$H$ be a subgroup of $Autf(M)$ of index strictly less than $2^{\aleph_{0}}$.
Then, there exists a finite subset $A\subset M$ such that $Autf_{A}%
(M)\subseteq H$.
\end{conjecture}

In \cite{lascar92}, this conjecture is proved for almost strongly minimal sets
(so, in particular for algebraically closed fields).

\section{The uncountable case}

\vskip-5mm \hspace{5mm}

We are now given a saturated structure $M$ of cardinality $\lambda>\aleph_{0}
$. The small index property has a natural generalization. If we assume that
$\lambda^{<\lambda}=\lambda$ (i.e. there is exactly $\lambda$ subsets of $M$
of cardinality less than $\lambda$) then any subgroup of $Aut(M)$ containing a subgroup of
the form $Aut_{A}(M)$ with $card(A)<\lambda$ has index at most $\lambda$. The
converse is true:

\begin{theorem}
\emph{\cite{lascarshelah}} Assume that $M$ is a saturated structure of
uncountable cardinality $\lambda=\lambda^{<\lambda}$, and let $H$ be a
subgroup of $Aut(M)$ of cardinality at most $\lambda$. Then, there exists a
subset $A$ of $M$ of cardinality less than $\lambda$ such that $Aut_{A}%
(M)\subseteq H$.
\end{theorem}

Here again, we may introduce a topological structure on $Aut(M)$: if $\mu$ is an
infinite cardinal, let $\mathcal{T}_{\mu}$ be the group topology on $Aut(M)$
for which a basis of open neighborhoods of the identity is
\[
\left\{  Aut_{A}(M)\;;\;A\subseteq M\text{ and }card(A)<\mu\right\}.
\]
To complete this definition, let $\mathcal{T}_{a}(M)$ be the group topology on
$Aut(M)$ for which a basis of open neighborhoods of the unit is
\[
\left\{  Autf_{A}(M)\;;\;A\subseteq M\text{ and }A\text{ finite}\right\}.
\]

The above theorem just says that the subgroups of $Aut(M)$ of index at most
$\lambda$ are exactly the open subgroups for $\mathcal{T}_{\lambda}(M)$, and
consequently, the topology $\mathcal{T}_{\lambda}(M)$ can be reconstructed
from the pure group structure. It is also clear that the open subsets for
$\mathcal{T}_{\lambda}$ are just the intersections of less than $\lambda$
$\mathcal{T}_{a}$-open subsets. So, if one knows $\mathcal{T}_{a}(M)$, one
knows $\mathcal{T}_{\lambda}(M)$.

With a few cardinality hypotheses, we can reconstruct one topological group
from another: (see \cite{lascar91} for details):

\begin{enumerate}
\item  Let $M$ and $M^{\prime}$ be two saturated models of the same theory.
Then we can reconstruct $(Aut(M^{\prime}),\mathcal{T}_{a})$ from
$(Aut(M),\mathcal{T}_{a})$.

\item  Let $M$ and $M^{\prime}$ be two models of the same theory, and assume
$card(M)=\lambda<\mu=card(M^{\prime})$. Then we can reconstruct $Aut(M^{\prime
})$ from $(Aut(M),\mathcal{T}_{\lambda})$ (and from $Aut(M)$ alone if
$\lambda=\lambda^{<\lambda})$. In fact we can reconstruct $(Aut(M^{\prime
}),\mathcal{T}_{\upsilon}(M^{\prime}))$ for every cardinal $\upsilon$,
$\lambda\leq\upsilon\leq\mu$.

\item  Let $M$ be a saturated structure of uncountable cardinality
$\lambda=\mu^{+}=2^{\mu}$ and assume that $T$ has a saturated model of
cardinality $\mu$. Then $(Aut(M),\mathcal{T}_{a})$ can be reconstructed from
$Aut(M)$.
\end{enumerate}

Let us give an example of a theorem which can be proved using the above facts: Assume GCH and
let $T_{1}$ and $T_{2}$ be two theories with saturated models $M_{1}$ and
$M_{2}$ of cardinality $\mu^{++}$. Assume that $Aut(M_{1})$ and $Aut(M_{2})$
are isomorphic. Then, for all cardinal $\lambda$, if $T_{1}$ has a saturated
model of cardinality $\lambda$, then $T_{2}$ has also a saturated model in
cardinality $\lambda$, and the automorphism groups of these two models are isomorphic.

\label{lastpage}

\end{document}